\documentclass{amsart}

\usepackage{amsthm,amssymb}
\usepackage[pagebackref,colorlinks,linkcolor=blue,citecolor=blue,urlcolor=blue,a4paper,hypertexnames=true]{hyperref}
\usepackage{amsrefs}
\usepackage[matrix, arrow]{xy}

\newcommand{\cU}{\mathcal U}
\newcommand{\N}{\mathbb N}
\newcommand{\Q}{\mathbb Q}
\newcommand{\Z}{\mathbb Z}
\newcommand{\C}{\mathbb C}
\newcommand{\calg}{\mathcal G}
\newcommand{\calk}{\mathcal K}
\newcommand{\calr}{\mathcal R}
\newcommand{\calh}{\mathcal H}
\newcommand{\Ext}{{\rm Ext}}

\newcommand{\Der}{Z^1}
\newcommand{\Inn}{B^1}

\newcommand{\ran}{{\rm ran}}
\newcommand{\res}{{\rm res}}

\theoremstyle{plain}
\newtheorem*{theorem*}{Theorem}

\newtheorem{theorem}{Theorem}[section]
\newtheorem{coro}[theorem]{Corollary}
\newtheorem{lemma}[theorem]{Lemma}
\newtheorem{propo}[theorem]{Proposition}

\theoremstyle{definition}

\newtheorem{definition}[theorem]{Definition}

\theoremstyle{remark}

\newtheorem{remark}[theorem]{Remark}
\newtheorem{example}[theorem]{Example}

\begin{document}
\title{Group cocycles and the ring of affiliated operators}

\author{Jesse Peterson}
\address{Jesse Peterson, University of California, Berkeley, U.S.A.}

\author{Andreas Thom}
\address{Andreas Thom, Universit\"at G\"ottingen, Germany.}

\begin{abstract}
In this article we study cocycles of discrete countable groups with values in $\ell^2 G$ and the ring of affiliated operators
$\cU G$. We clarify properties of the first cohomology of a group $G$ with coefficients in $\ell^2 G$ 
and answer several questions from \cite{ctv}.
Moreover, we obtain strong results about the existence of free subgroups and the subgroup structure, 
provided the group has a positive first $\ell^2$-Betti number. We give numerous applications and examples
of groups which satisfy our assumptions.
\end{abstract}

\maketitle


\section{Introduction} \label{intro}

Let $G$ be a discrete countable group and let $M$ be a $G$-module. A \emph{cocycle} $c\colon G \to M$ is a map which
satisfies $$c(gh) = g\cdot c(h) + c(g).$$
It is called \emph{inner}, if there exists $\xi \in M$, such that $c(g) = (g-1)\xi$. We denote by $Z^1(G;M)$ the space
of cocycles, by $B^1(G;M)$ the subspace of inner cocycles and by $H^1(G;M)$ the first cohomology of the group $G$ with coefficients in $M$, i.e.\ the quotient of $Z^1(G;M)$ by $B^1(G;M)$.
Many properties of $G$ can be phrased in terms of cocycles and cohomology; and with this article we support
the view that certain $G$-modules of functional analytic nature turn out to be particularly useful in the study of infinite groups. The closer study of the special case $M = \ell^2G$ is usually called the theory of $\ell^2$-invariants of groups. 

Here, we denote by $\ell^2 G$ the Hilbert space with basis $G$, and by $B(\ell^2G)$ the Banach space
of bounded linear endomorphisms of $\ell^2G$.
The left and right translation of $G$ on itself extend to two commuting unitary representations:
$$\lambda,\rho \colon G \to U(\ell^2 G) = \{u \in B(\ell^2G) \mid uu^* =u^*u=1\},$$
and endow $\ell^2G$ with a left and a right $G$-module structure. 
It is well-known that the generated von Neumann algebras $LG = \lambda(G)''$ and $RG = \rho(G)''$
are commutants of each other. Here $S' = \{t \in B(\ell^2G) \mid st=ts, \ \forall s\in S\}$ denotes 
the commutant of the set $S$. $LG$ is called the group von Neumann algebra of $G$. We frequently identify
$RG$ with $LG$ and consider $\ell^2G$ as a $LG$-bimodule.

The theory of $\ell^2$-invariants was started in the seminal work of M.\ Atiyah in \cite{atiyah} and developed further by 
J.\ Dodziuk, see \cite{dodziuk}, J.\ Cheeger and M.\ Gromov in \cite{cheegergromov}. Among many others, major contributions were obtained by W.\ L\"uck, see \cite{lueck}, and D.\ Gaboriau, see \cite{gab2}.

In our study, the ring $\cU G$ of densely defined, closed
operators on $\ell^2G$, which are affiliated with $LG$, is of major importance. For details about its
definition consult \cite{tak2}*{Chapter IX}. We naturally have the following chain of inclusions of $G$-modules:
$LG \subset \ell^2 G \subset \cU G$, which induce maps on cocycles and cohomology.

\vspace{0.2cm}

The first $\ell^2$-Betti number $\beta_1^{(2)}(G)$ is defined to be a certain dimension of either 
$H_1(G,\ell^2G)$ or $H^1(G,\ell^2 G)$, see Section \ref{l2coh}. 
It turns out to be useful to study $Z^1(G;\ell^2G)$ through its map to $Z^1(G;\cU G)$. In the case where the group $G$ in non-amenable we show that the first $\ell^2$-Betti number vanishes if and 
only if $H^1(G, \ell^2G) = 0$ which was previously shown for finitely generated groups in \cite{bv}.  In case the group $G$ is amenable,
$Z^1(G;\cU G)=B^1(G;\cU G)$, and we can show that each element in $c  \in Z^1(G;\ell^2 G)$ is either bounded
or proper on $G$, depending on whether the vector $\xi \in \cU G$, for which $c(g)=(g-1)\xi$, is in 
$\ell^2G$ or not.
\begin{theorem*} (see Theorem \ref{propb}) Let $G$ be an countable and discrete group which is amenable.
Every $1$-cocycle with values in $\ell^2 G$ is either bounded or proper.
\end{theorem*}
Moreover, the existence of co-cycles which are neither proper nor bounded is proved for
non-amenable $G$, with the necessary condition of non-vanishing first $\ell^2$-Betti number, and 
provided there exists an infinite amenable subgroup of $G$. 

\vspace{0.2cm}

Providing a large of examples of groups with a positive first $\ell^2$-Betti number, we prove:

\begin{theorem*} (see Theorem \ref{posbetti}) Let $G$ be an infinite 
countable discrete group. Assume that 
$$G = \langle g_1, \dots , g_n \mid  r_1^{w_1},\dots,r_k^{w_k} \rangle,$$
for elements $r_1,\dots,r_k \in {\mathbb F}_n = \langle g_1,\dots,g_n \rangle$ and positive integers $w_1,\dots,w_k$.
We assume that the presentation is irredundant in the sense that $r_i^l \neq e \in G$, for $1<l < w_i$ and $1 \leq i \leq k$.
Then, the following inequality holds:
$$\beta_{1}^{(2)}(G) \geq n-1 - \sum_{j=1}^k \frac{1}{w_j}. $$
\end{theorem*}

Denis Osin has used this Theorem in \cite{MR2552302} to construct $n$-generated torsion groups with first $\ell^2$-Betti number greater than $n-1-\varepsilon$.

\vspace{0.2cm}

We study the existence
of free subgroups in torsionfree discrete groups. It is well-known that non-amenability is not sufficient to ensure
the existence of a free subgroup. We show that a non-vanishing first $\ell^2$-Betti number is sufficient, provided
$G$ is torsionfree and satisfies a weak form of Atiyah's Conjecture, see Section \ref{freegroups}. More precisely,

\begin{theorem*} (see Theorem \ref{st1})
Let $G$ be a torsionfree discrete countable group.
There exists a family of subgroups $\{ G_i \mid i \in I \}$, such that
\begin{enumerate}
\item[(i)] We can write $G$ as the disjoint union:
$$G = \{e\} \cup \bigcup_{i \in I} \dot G_i.$$
\item[(ii)] The groups $G_i$ are mal-normal in $G$, for $i \in I$.
\item[(iii)] If $G$ satisfies a weak form of the Atiyah Conjecture, then $G_i$ is free from $G_j$, for $i \neq j$.
\item[(iv)] $\beta_1^{(2)}(G_i)=0$, for all $ i \in I$.
\end{enumerate}
\end{theorem*}

Moreover, we
obtain a strong structure theorem for such groups. The techniques
allow to generalize a recent result of J.\ Wilson, see \cite{MR2011971}. This also gives a new estimate on the
exponential growth rate in terms of the first $\ell^2$-Betti number; proving a generalized form of 
Conjecture $5.14$ of Gromov from \cite{MR1699320}.

\vspace{0.2cm}

D.\ Gaboriau proved in \cite{gab2} that
an infinite, normal, infinite index subgroup $H$ of a group $G$ with positive first $\ell^2$-Betti 
number cannot have a finite first $\ell^2$-Betti number, and in particular cannot be finitely generated. 
Assuming infinite index, we extend this result to subgroups $H$, for which $H \cap H^g$ is infinite, for 
all $g \in G$. (See Section \ref{notionsofnormality} for the defintion of $wq$-normality and $ws$-normality.) In particular, it applies to all subgroups which contain an infinite normal subgroup. 
This covers classical results by Karrass-Solitar \cite{MR0086813}, Griffiths \cite{MR0153832} and Baumslag \cite{MR0199247}, as well as more recent results
by Bridson-Howie \cite{MR2262789} and Kapovich \cite{MR1859278}.
More precisely, we prove:

\begin{theorem*}(see Theorems \ref{main1} and \ref{main2})
Let $G$ be a countable discrete group and suppose $H$ is an infinite subgroup.
\begin{enumerate}
\item If $H$ $wq$-normal subgroup, then $\beta_1^{(2)}(H) \geq \beta_1^{(2)}(G)$.
\item If $ws$-normal and infinite index, and $\beta_1^{(2)}(H) < \infty$, then $\beta_1^{(2)}(G) =0$.
\end{enumerate}
\end{theorem*}

Among the corollaries, we prove that if 
$H,K$ are infinite, finitely generated subgroups of $G$, so that $H \cap K$ is of finite index in $H$ and $K$, 
then: the index of $H \cap K$ in $\langle H,K \rangle$ is finite, if the first $\ell^2$-Betti number of $\langle H,K \rangle$ is non-zero, see Theorem \ref{finiteindex}.
Many of these results are well-known for free groups and were proved by several authors in various other 
cases, which are mostly covered by our result. The particular case of limit groups was studied in \cite{MR1859278}. 
Our proof is using concrete computations
with cocycles with values in $\cU G$ and results from ergodic theory.
It was D.\ Gaboriau in his groundbreaking work \cite{gab2}, 
who was the first to use ergodic theory to obtain striking results in the theory of $\ell^2$-invariants with
applications to infinite groups.
\vspace{0.2cm}

The article is organized as follows:

\vspace{0.1cm}

Section \ref{intro} is the Introduction. 
In Section \ref{l2coh} we recall the program of W.\ L\"uck and introduce
the algebra $\cU G$ of densely defined closed operators, which are affiliated with the group von Neumann algebra 
$LG$. Several algebraic properties of $\cU G$ are recalled and their implications are clarified. In Theorem
\ref{cohom}, we show that L\"uck's generalized dimension of the first cohomology with coefficients in either
$LG, \ell^2G$ or $\cU G$ coincides with the first $\ell^2$-Betti number. Moreover, for $H \subset G$, we show that 
$H^1(H;\cU G)=0$, whenever the first $\ell^2$-Betti number of $H$ vanishes.
This will turn out to be very useful in algebraic computations. 

Using these results we prove
that a $\ell^2$-cocycle on an amenable group is either bounded or proper. This is Conjecture $2$ from \cite{ctv}.
Theorem \ref{neitherbounded} clarifies the existence of co-cycles which are neither bounded nor proper for
general groups (admitting an infinite amenable subgroup).

\vspace{.1cm}

In Section \ref{examples} we give examples of groups with non-vanishing first $\ell^2$-Betti number.
We give lower bounds on the first $\ell^2$-Betti numbers for amalgamated free products, HNN-extensions
and various other more elaborate constructions. We hope that this section 
provides some useful tools, for example Theorem \ref{posbetti}, to 
estimate the first $\ell^2$-Betti in some interesting cases.

\vspace{.1cm}

In Section \ref{freegroups}, we examine torsionfree groups which satisfy a weak form of the Atiyah Conjecture.
It turns out that a positive first $\ell^2$-Betti number has strong implications on the structure of such groups. In
Theorem \ref{st1}, we show that any such group decomposes \emph{as a pointed set} into malnormal, mutually
free subgroups with vanishing first $\ell^2$-Betti number. This result is used to prove that in this case 
the reduced group $C^*$-algebra is simple with a unique trace state. Moreover, the techniques imply a Freiheitssatz, 
see Corollary \ref{manyfreesubgroups},
and give estimates about the exponential growth rate. In particular, assuming a weak form of Atiyah's Conjecture, 
we obtain a new proof of Conjecture $5.14$ of M.\ Gromov in \cite{MR1699320} about the exponential growth rate of a finitely presented group with fewer relations than generators.

\vspace{.1cm}

Section \ref{proofsec} contains the main results of this article. We introduce various notions of normality, see
the definition in Section \ref{notionsofnormality}, and study the existence of infinite subgroups of a group $G$ 
with infinite index, sharing one of the normality properties, see Theorems \ref{main1} and \ref{main2}. In particular, if  $\beta_1^{(2)}(G) \neq 0$, we can exclude the existence of a finitely generated subgroup of infinite index, which contains an infinite normal subgroup. Several other corollaries can be found in this section, whereas various other applications of the main theorems are contained in Section \ref{misc}.
The proof of Theorem \ref{main2} relies on discrete measured groupoids and ergodic theory. 

\vspace{.1cm}

The necessary results from ergodic theory and the theory of discrete measured groupoids 
are collected in Section \ref{ergodic}, where we extend some of our results from Section \ref{proofsec}
to discrete measured groupoids.

\vspace{.1cm}
 
In Section \ref{misc} we collect results about various classes of groups. In particular, we study boundedly 
generated groups, certain groups which are generated by a family of subgroups, limit groups, groups which are measure equivalent to free groups and so-called
powerabsorbing subgroups. We are able to reprove and generalize several results from the literature. Most notably, we prove
Proposition \ref{finiteindex}, which is a generalization of Theorem $C$ of I.\ Kapovich in \cite{MR1859278}.

\section{$\ell^2$-cohomology, cocycles and Betti numbers} \label{l2coh}

The computations of $\ell^2$-homology have been algebraized through the seminal work of W.\ L\"uck, which is
summarized and explained in detail in his nice compendium \cite{lueck}. The basic observation is that 
through a dimension function, which is defined for all modules over the group von Neumann algebra, entirely algebraic objects give rise to numerical invariants. One of the main results is the following equality:
\begin{equation}
\beta_n^{(2)}(G) = \dim_{LG} H_n(G;LG),
\end{equation}
where $G$ is a countable discrete group and $\beta_n^{(2)}(G)$ denotes the $n$-th $\ell^2$-Betti number in the sense
of M.\ Atiyah and as generalized by J.\ Cheeger and M.\ Gromov, see \cites{atiyah, cheegergromov}. We are freely
using $\dim_{LG}$, L\"uck's dimension function, which is defined for all $LG$-modules. We are frequently using that for an
extension $$0 \to M' \to M \to M'' \to 0$$
of $LG$-modules, the following formula holds:
$$\dim_{LG} M = \dim_{LG} M' + \dim_{LG} M''.$$
For more details about the definition of the dimension function and its properties consult \cite{lueck}.

In our study of $\ell^2$-cohomology, the ring of closed and densely defined operators on $\ell^2G$, which are affiliated
with $LG$, plays a prominent role. This ring is denoted by $\cU G$. Those rings were the motivating examples for
J.\ von Neumann to study rings with a remarkably strong algebraic property, which was later named \emph{von Neumann regularity}. A ring $R$ is said to be von Neumann regular, if for each $a \in R$, there exists $b \in R$, such that $aba=a$.
Alternatively, von Neumann regular rings are precisely those, for which all modules are flat. Recall, a module $M$ over a ring $R$ is called \emph{flat}, if the functor $? \otimes_R M$ is exact.

In the process of algebraization of $\ell^2$-homology, it was P.\ Linnell in \cite{linnell}, who reintroduced the ring of affiliated operators and studied its nice algebraic properties. Some ring theoretic properties were studied towards applications
to $\ell^2$-invariants and $K$-theory in \cites{reich, thoml2}.

To our knowledge, this following lemma was first observed by S.K. Berberian (see \cite{berberian}) in 
the context of finite von Neumann algebras and shortly afterwards by K.R.\ Goodearl in \cite{goodearl} in 
the more general context of metrically complete von Neumann regular rings. 

\begin{lemma} \label{goode} Let $(M,\tau)$ be a finite tracial von Neumann algebra. The ring $\cU G$ of operators affiliated with $M$ is self-injective. \end{lemma}

Recall, a ring $R$ is called \emph{self-injective}, if the functor $\hom_R(?,R)$ is exact. (Note that $R \cong R^{\rm op}$ for all our rings, so that we do not need to talk about \emph{left} self-injectivity etc.)
Although, Lemma \ref{goode} has been around for more than $20$ years, its consequences for the computation of $\ell^2$-invariants have not been fully exploited. There are indications, that the context of metrically complete modules over metrically complete rings is indeed a useful context to study $\ell^2$-invariants. Indeed, in \cite{thomrank}, the second author gave a conceptual and short proof of D.\ Gaboriau's result about invariance of $\ell^2$-Betti numbers under orbit equivalence. Moreover in \cite{sauerthom}, using essential properties of the category of metrically complete modules, R.\ Sauer and the second author constructed a Hochschild-Serre spectral sequence for extensions of discrete measurable groupoids.

In this note, we are basically interested in the first $\ell^2$-Betti number of a countable discrete group.
Right from the beginning of the study of $\ell^2$-homology and $\ell^2$-cohomology of discrete groups, 
it was observed that those are only dual to each other under some finiteness assumptions on the group, i.e.\ the group was
assumed to have a finite classifying space. However, using the self-injectivity of $\cU G$, it
is obvious that always:
\begin{equation}
\hom_{\cU G}\left(H_n(G; \cU G),\cU G\right) \cong H^n(G; \cU G).
\end{equation}

Moreover, $LG \subset \cU G$ is a flat ring extension and $? \otimes_{LG} \cU G$ preserves the dimension, see
\cites{reich, thomrank}. Hence
$$\dim_{LG} H_n(G,LG) = \dim_{LG} H_n(G,LG) \otimes_{LG} \cU G = \dim_{LG} H_n(G,\cU G).$$
Note that,  by Corollary $3.4$ in \cite{thoml2}, also dualizing a $\cU G$-module preserves its dimension.
We conclude that
\begin{equation} \beta_n^{(2)}(G) = \dim_{LG} H^n(G; \cU G).\end{equation}
The computations with this cohomology group simplify the picture drastically since they have the nice
property that they vanish unless their dimension is non-zero. This follows from Corollary $3.3$ in \cite{thoml2}.
We will use this fact frequently.
\subsection{A cocycle description}

It is well-known that the first group cohomology with coefficients in a module $M$ 
can be computed as the vector space of $M$-valued $1$-cocycles on $G$ modulo inner cocycles. A $1$-cocycle with values in the $G$-module $M$ is a map
$$ c\colon G \to M, \quad \mbox{with} \quad c(gh) = gc(h) + c(g), \quad \forall g,h \in G.$$
It is called inner, if there exists $\xi \in M$, such that $c(g) = (g-1) \xi$, for all $g \in G$.
We denote the space of $M$-valued $1$-cocycles by $\Der(G;M)$ and the space of inner cocycles
by $\Inn(G;M)$. There is an exact sequence
$$0 \to \Inn(G;M) \to \Der(G;M) \to H^1(G;M) \to 0.$$

Our first theorem gives an identification of dimensions of cohomology groups, where the coefficients
vary among the canonical choices $LG,\ell^2G$ and $\cU G$.

\begin{theorem} \label{cohom}
Let $G$ be a countable discrete group.
\[\beta_k^{(2)}(G) = \dim_{LG} H^k(G,\cU G) = \dim_{LG} H^k(G,\ell^2 G) = \dim_{LG} H^k(G,L G).\]
Moreover, if $\beta_k^{(2)}(G) = 0$ for some $k$, then $H^k(G,\cU G)=0$.
\end{theorem}
\begin{proof} We give the proof only for $k=1$, since we are mainly concerned with the first $\ell^2$-cohomology. A similar argument can be found for $k \neq 1$. There exists a commutative diagram with exact rows as follows:
$$
\xymatrix{ 
0 \ar[r] & \Inn(G;LG) \ar[r] \ar[d] & \Der(G;LG) \ar[r] \ar[d] & H^1(G;LG) \ar[r] \ar[d]& 0 \\
0 \ar[r] & \Inn(G;\ell^2G) \ar[r] \ar[d] & \Der(G;\ell^2G) \ar[r] \ar[d] & H^1(G;\ell^2G) \ar[r] \ar[d]& 0 \\
0 \ar[r] & \Inn(G;\cU G) \ar[r]  & \Der(G;\cU G) \ar[r]  & H^1(G;\cU G) \ar[r] & 0. \\
}
$$
Recall, if $M_1 \subset M_2$ is an inclusion of $LG$ module, then it is called \emph{rank dense}, if for every
$\xi \in M_2$, there exists an increasing sequence of projections $p_n \uparrow 1$, such that $\xi p_n \in M_1$, for
all $n \in \N$. It was shown in \cite{thomrank}, that a rank dense inclusion is a dimension isomorphism.

If $G$ is infinite, the left column identifies with the inclusions $LG \subset \ell^2 G \subset \cU G$, which are
well-known to be dimension isomorphisms since $LG$ is rank dense in $\cU G$.

The column in the middle also consists of inclusions and we claim that the images are rank dense as well. 
Indeed, every $1$-cocycle with values in $\cU G$ can be cut by a projection 
of trace bigger than $1 - \varepsilon$ to take values in $LG$. 
For each $g \in G$, $c(g) \in \cU G$ and we find a projection $p_g$ of 
trace bigger $1-\varepsilon_g$, so that $c(g)p_g \in LG$. 
Taking the infimum over all $p_g$, we obtain a projection $p$ of trace bigger than
$1 -\sum_{g \in G} \varepsilon_g$. Hence, choosing a suitable sequence $\varepsilon_g$ 
proves the claim. 

The vanishing of $H^1(G;\cU G)$ in case of vanishing first $\ell^2$-Betti number follows since it is the dual of
the $\cU G$-module $H_1(G;\cU G)$. It was shown in \cite{thoml2}, that the dual is zero if and only if the dimension
is zero. This finishes the proof.
\end{proof}

\begin{remark}
Let $H \subset G$ be a subgroup. 
It follows from standard computations that $$\beta_1^{(2)}(H) = \dim_{LG}H^1(H,\cU G)$$ and $H^1(H,\cU G) =0$,
if and only if the first $\ell^2$-Betti number of $H$ vanishes.
\end{remark}

In \cite{bv} it is shown that for a finitely generated non-amenable discrete group, the first $\ell^2$-Betti number vanishes if
and only if the first cohomology group with values in the left regular representation vanishes, (see also Corollary 3.2 in
\cite{mv}). We will now show that we may drop the assumption that the group is finitely generated.

\begin{coro} Let $G$ be a non-amenable countable discrete group,
then $\beta_1^{(2)}(G) = 0$ if and only if $H^1(G, \ell^2G) = 0$.
\end{coro}
\begin{proof} First let us suppose that $H^1(G, \ell^2G) = 0$. Let $c: G \to \cU G$ be a $1$-cocycle, we must
show that there is an affiliated operator $\xi \in \cU G$ such that $c(g) = (g - 1)\xi$, for each $g \in G$.
Given $\varepsilon > 0$, an affiliated operator $\eta$ and a projection $p \in RG$ we may find a projection
$q \in RG$, $q \leq p$ such that $q\eta \in \ell^2 G$ and $\tau(p - q) < \varepsilon$.  From this fact
we may construct a partition of unity $\{ p_n \}_{n \in \N}$ in $RG$ such that
$$p_n c(g) \in \ell^2 G,\quad \forall g \in G, n \in \N.$$

Since $H^1(G, \ell^2G) = 0$ we conclude that there exist $\xi_n \in \ell^2 G$
such that $$p_n c(g) = (g - 1)\xi_n, \quad \forall g \in G, n \in \N.$$
Moreover we may assume $p_n \xi_n = \xi_n$, $\forall n \in \N$ so
that $\xi = \sum_{n \in \N} \xi_n \in \cU G$ is well defined and has the desired properties.

Next let us suppose that $\beta_1^{(2)}(G) = 0$.  From Theorem \ref{cohom},
we conclude that if $c:G \to \ell^2G$ is a $1$-cocycle then there exists
and affiliated operator $\xi \in \cU G$ such that $c(g) = (g - 1)\xi$, $\forall g \in G$.
Let $\{ p_n \}_{n \in \N}$ be a sequence of projections in $RG$ which increase to $1$ such that
$p_n \xi \in \ell^2G$, for each $n \in \N$. Then since $RG$ acts normally on $\ell^2G$ we conclude that
$$\lim_{n \to \infty} \| (1 - p_n) c(g) \| = 0,\quad \forall g \in G.$$  Hence $c$ is approximately inner
which shows that $\overline{H^1}(G, \ell^2G) = 0$.  As $\ell^2G$ does not weakly contain
the trivial representation it is a well known result that we must also have that
$H^1(G, \ell^2G) = 0$.
\end{proof}

\subsection{Dichotomy of $\ell^2$-cocycles on amenable groups}

The following theorem is an affirmative answer to Conjecture $2$ in \cite{ctv}.

\begin{theorem} \label{propb} Let $G$ be an countable and discrete group which is amenable.
Every $1$-cocycle with values in $\ell^2 G$ is either bounded or proper.
\end{theorem}
\begin{proof} Let $c: G \to \ell^2 G$ be a $1$-cocycle. We need to show that either $\sup_{g \in G} \|c(g)\|_2$ is finite or
$\{\|c(g_n)\|_2, n \in \N\}$ is unbounded for every sequence $\{g_n\}_{n \in \N}$ in $G$ that goes to infinity.

It is well-known, that $\beta_1^{(2)}(G) = 0$ if $G$ is amenable. We conclude from Theorem \ref{cohom},
that there exists an affiliated operator $\xi \in \cU G$, such that $$c(g) = (g-1)\xi,$$ forall $g \in G$.
Let $\{p_m\}_{m \in \N}$ be a partition of unity of $RG$, which has the property that $p_m \xi \in \ell^2 G$, for all $m \in \N$.

Since the left-regular representation is mixing, for every sequence $\{g_n\}_{n \in \N}$ that goes to infinity
\[\lim_{n \to\infty} \|(g_n-1)p_m \xi\|^2_2 = 2 \|p_m\xi\|^2_2.\]

We consider two cases depending whether $\sum_{m \in \N} \|p_m\xi\|_2^2$ is finite or not.
(a) If it is finite, then $\xi \in \ell^2 G$ and the cocycle will be bounded.
(b) If it is infinite we aim to show that the cocycle is proper. Let $\{g_n\}_{n \in \N}$ be a sequence in $G$ that tends to infinity.
Given $C >0$, there exists $k \in \N$, such that
\[\sum_{m=1}^k 2\|p_m\xi\|^2_2 \geq C + 1,\]
and there exists some $l \in \N$, such that $$\|(g_j -1)p_m\xi\|_2^2 \geq 2\|p_m\xi\|_2^2 - k^{-1}$$ for all $1 \leq m \leq k$ and all $j \geq l$.
This implies that for every $j \geq l$ we have
\[
\|(g_j-1)\xi\|^2_2 \geq \sum_{m=1}^k \|(g_j-1)p_m\xi\|^2_2 \geq \sum_{m=1}^k 2 \|p_m \xi\|_2^2 - k^{-1} \geq C.
\]
This finishes the proof.
\end{proof}

We remark that it was previously shown in \cite{mv} that if $G$ is a countable discrete group and $c:G \to \ell^2G$ is
an unbounded $1$-cocycle then $c$ is also unbounded on any infinite subgroup of $G$.  The above theorem states that this is true for infinite subsets as well.
There is a partial converse to the preceding result.
\begin{theorem} \label{neitherbounded}
Let $G$ be a group with $\beta^{(2)}_1(G) \neq 0$ and assume that there exists an infinite amenable sub-group.
There exists a $1$-cocycle with values in $\ell^2 G$ on $G$ which is neither bounded nor proper.
\end{theorem}

\begin{proof}
Let $H \subset G$ be an infinite amenable subgroup.  Since $\beta^{(2)}_1(H) = 0 < \beta^{(2)}_1(G)$ the restriction map
$H^1(G,\ell^2 G) \to H^1(H, \ell^2 G)$ cannot be injective.  Hence there is an unbounded $\ell^2$-cocycle on $G$ which is bounded on $H$.
\end{proof}

\begin{remark}
Note that for a non-amenable group with vanishing first $\ell^2$-Betti number, all $\ell^2$-cocycles are automatically bounded.
\end{remark}

The following corollary answers Question $1$ in \cite{ctv} negatively.

\begin{coro} If $G = H_1 \ast H_2$ with $H_j$ non-trivial, $j=1,2$ and $H_1 \neq \Z/2\Z$, then there exists
an $\ell^2$-cocycle, which is neither proper nor bounded. In particular, this is the case for $PSL_2(\Z) = \Z/2\Z \ast \Z/3\Z$.
\end{coro}

\section{Examples of groups with positive first $\ell^2$-Betti number} \label{examples}

This class of groups (or rather its complement) was also studied by B.\ Bekka and A.\ Valette in \cite{bv}.
Among the classical examples, there are free groups, surface groups and groups containing such groups with finite index. W.\ Dicks and P.\ Linnell (see \cite{MR2285740}) computed that the first $\ell^2$-Betti number of a $n$-generated one-relator group is $n-2$.
It was shown in \cite{MR2250857}, that the class of finitely generated groups with first $\ell^2$-Betti number greater or equal
than $\varepsilon$ is closed in Grigorchuk's space of marked groups. This implies in particular, that limit groups have positive
first $\ell^2$-Betti number, see Section \ref{limitgroups}.

D.\ Gaboriau showed in \cite{gab2} that the non-vanishing of the $n$-th $\ell^2$-Betti number does only depend on
the group up to measure equivalence. This provides a class of examples which we study more closely in Section 
\ref{measeq}.

\vspace{0.2cm}

Throughout this section, we are using the convention that $|G|^{-1}=0$, if $G$ is infinite.

\subsection{Amalgamated free products}
In the case of amalgameted free products, we state the following well-known result:
\begin{propo} Let $G$ be a discrete countable group.
If $G$ is an amalgamated free product $G = A \ast_C B$, then 
\begin{equation} \label{amalga}
\beta_1^{(2)}(G) \geq \left(\beta_1^{(2)}(A) - \frac1{|A|} \right)+ \left( \beta_1^{(2)}(B) - \frac1{|B|} \right) 
- \left( \beta_1^{(2)}(C) - \frac1{|C|} \right).
\end{equation}
\end{propo}
\begin{proof}
This follows from the long exact sequence of homology with coefficients in $\cU G$, which is associated to
an amalgamated free product (see \cite{brown}*{Ch.\ VII.9}) and an easy dimension count.
\end{proof}

Note, if $G$ acts on a tree, similar estimates can be found in terms of the order and the first $\ell^2$-Betti numbers
of the stabilizer groups. In particular, if $G = A\ast_B$ is an HNN-extension, then:
$$\beta_1^{(2)}(G) \geq \left(\beta_1^{(2)}(A) - \frac1{|A|} \right)- \left( \beta_1^{(2)}(B) - \frac1{|B|} \right) .$$
Of course, in special cases, much more can be said.  For example, 
from Apendix A of \cite{bdj} we see that the inequality in
\ref{amalga} is actually an equality in the case where $\beta_1^{(2)}(C) = 0$.

\subsection{Triangle groups and related constructions}

We now want to provide more non-trivial examples of groups which have a positive first $\ell^2$-Betti number.
Our main result in this direction is the following theorem:

\begin{theorem} \label{posbetti} Let $G$ be an infinite 
countable discrete group. Assume that there exist subgroups $G_1,\dots,G_n$,
such that
$$G = \langle G_1, \dots , G_n \mid  r_1^{w_1},\dots,r_k^{w_k} \rangle,$$
for elements $r_1,\dots,r_k \in G_1 \ast \cdots \ast G_n$ and positive integers $w_1,\dots,w_k$.
We assume that the presentation is irredundant in the sense that $r_i^l \neq e \in G$, for $1<l < w_i$ and $1 \leq i \leq k$.
Then, the following inequality holds:
$$\beta_{1}^{(2)}(G) \geq n-1 + \sum_{i=1}^n \left( \beta_1^{(2)}(G_i) - \frac{1}{|G_i|}\right) - \sum_{j=1}^k \frac{1}{w_j}. $$
\end{theorem}
\begin{proof}
There is an exact sequence
$$ 0 \to Z^1(G;\cU G) \stackrel{p}{\to} Z^1(G_1 \ast \cdots \ast G_n;\cU G) \stackrel{q}{\to} \oplus_{j=1}^k \cU G,$$ 
where $p$ is given by the composition
$$Z^1(G;\cU G) \to \oplus_{i=1}^n Z^1(G_i;\cU G) \cong Z^1(G_1 \ast \cdots \ast G_n;\cU G),$$
and $q$ is given by the sum of the evaluation maps at $r_j^{w_j}$.
Indeed, exactness at $Z^1(G;\cU G)$ is clear and it remains to prove exactness in the middle. Again, it is obvious
that the composition is zero so that we only have to show that any element in the kernel of the evaluation maps defines
a cocycle on $G$. If a cocycle on $G_1 \ast \cdots \ast G_n$ vanishes on a relator $r_i^{w_i}$, then it also vanishes
on any of its conjugates, since:
$$c(gr_i^{w_i} g^{-1}) = (1-gr_i^{w_i} g^{-1})c(g) + g c(r_i^{w_i}) = (1-gr_i^{w_i} g^{-1})c(g) =0.$$
Here, we are using that $gr_i^{w_i} g^{-1}$ acts trivially on $\cU G$.

Clearly, $$c\left(r_j^{w_j}\right) = \sum_{l=1}^{w_j} r_j^l \cdot c(r_j).$$ Thus the dimension of the image of the evaluation map 
at $r_j^{w_j}$ is less than $1/w_j$, since $1/{w_j} \cdot \sum_{l=1}^{w_j} r_j^l$ is a projection of trace $1/{w_j}$. Here, we are
using that the presentation is irredundant.
Hence, the dimension of the image of $q$ is less than $\sum_{j=1}^k 1/{w_j}$. 
The claim follows by noting that 
$$\dim_{LG} Z^1(G_i;\cU G) = \beta_{1}^{(2)}(G_i) - \frac1{|G_i|} +1.$$
\end{proof}

The theorem covers generalized triangle groups and so-called generalized tetrahedron groups. Let us spell out what
the theorem says in the case of generalized triangle groups. Let us first recall the definition.

\begin{definition} A group $G$ is called a \emph{generalized triangle group}, if it admits a representation
$$G = \langle a,b \mid a^p=b^q = w(a,b)^r \rangle,$$
where $w(a,b)$ is a cyclically reduced word of length at least $2$ in $C_p \ast C_q$.
We call $$\kappa(G) = \frac1p + \frac1q + \frac1r -1$$
the curvature of $G$.
\end{definition}

The following result is an immediate consequence of Theorem \ref{posbetti}.

\begin{coro}
Let $G$ be a triangle group. 
The following inequality holds:
$$\beta_1^{(2)}(G) \geq  - \kappa(G).$$ 
In particular, if $G$ is negatively curved, i.e.\ $\kappa(G)<0$, then $\beta_1^{(2)}(G) \neq 0$.
\end{coro}

\subsection{The relation module}
Another result which also estimates the first $\ell^2$-Betti number in terms of more algebraic data is
given by the following theorem:

\begin{propo} Let $G$ be a finitely generated group with a presentation $0 \to R \to F \to G$
with $F$ free of rank $n$.
Then the following inequality holds:
$$\beta_1^{(2)}(G) \geq n-1 - \dim_{LG} \left( R_{ab} \otimes_{\Z G} \cU G \right),$$
where $R_{ab} = R/[R,R]$ is the relation $G$-module (induced by the conjugation action of $G$ on $R$).
\end{propo}
\begin{proof}
Underlying the computation in the proof of Theorem \ref{posbetti}, there is a Lyndon-Serre spectral sequence with
a low degree exact sequence. Let $0 \to R \to F \to G \to 0$ be a presentation of the group $G$ with  $F$ 
free of rank $n$. Then,
$$0 \to H^1(G; \cU G) \to H^1(F; \cU G) \to H^1(R; \cU G)^{G} \to H^2(G;\cU G) \to 0$$
is an exact sequence. Note that $\cU G$ is a trivial $R$-module, so that
$$H^1(R;\cU G) = \hom_{\Z}(R_{ab} \otimes_{\Z} \C,\cU G),$$
where $G$ acts diagonally with the conjugation action on $R$ and on the left on $\cU G$.
Hence,
$$H^1(R;\cU G)^G = \hom_{\Z G}(R_{ab},\cU G).$$
Writing everything out, we get:
$$0 \to H^1(G; \cU G) \to \cU G^{ \oplus n-1} \to \hom_{\Z G}(R_{ab},\cU G) \to H^2(G;\cU G) \to 0.$$
Taking dimensions, this implies the claim.
\end{proof}

\section{Free subgroups} \label{freesubgroups}

\label{freegroups}
\subsection{Restriction maps and free subgroups}
Throughout this section, we are assuming that $G$ is a torsionfree discrete countable group and most of the time 
also that it satisfies the following condition:
\begin{equation*}
(\star) \qquad \mbox{Every non-trivial element of $\Z G$ acts without kernel on $\ell^2 G$.}
\end{equation*} 
This condition is satisfied if $G$ satisfies 
the Atiyah conjecture but is \emph{a priori} weaker. Recall, the Atiyah Conjecture for torsionfree groups
predicts the existence of a skew-field $\Z G \subset K \subset \cU G$. 
Note that the Atiyah Conjecture was established for a large class of
torsionfree groups, see the results and references in \cite{MR1990479}. In particular, condition $(\star)$ is known to hold 
for all right orderable groups and all residually torsionfree elementary amenable groups. The subset of Grigorchuk's space
of marked groups for which the conjecture holds is closed and hence contains for example all limit groups.
No counterexample is known.

\vspace{.1cm}

Let $G$ be a discrete group, we use the notation $\dot G$ to denote the set $G\setminus \{e\}$. In our computations we
are exploiting the basic fact that $1-g$ is invertible as an affiliated operator if $g$ is not torsion. Of course, this observation
does not rely on condition $(\star)$.
The main result here is the following theorem.

\begin{theorem} \label{st1}
Let $G$ be a torsionfree discrete countable group.
There exists a family of subgroups $\{ G_i \mid i \in I \}$, such that
\begin{enumerate}
\item[(i)] We can write $G$ as the disjoint union:
$$G = \{e\} \cup \bigcup_{i \in I} \dot G_i.$$
\item[(ii)] The groups $G_i$ are mal-normal in $G$, for $i \in I$.
\item[(iii)] If $G$ satisfies condition $(\star)$, then $G_i$ is free from $G_j$, for $i \neq j$.
\item[(iv)] $\beta_1^{(2)}(G_i)=0$, for all $ i \in I$.
\end{enumerate}
\end{theorem}
\begin{proof}
We partition $\dot G$ according to the following equivalence relation
$$g \sim h \quad \Leftrightarrow \quad c(g) =0 \quad \mbox{if and only if} \quad c(h) =0, \quad \forall c \in 
Z^1(G; \cU G).$$
First of all, one direction is sufficient to imply the \emph{if and only if} in the definition of the equivalence
relation. Indeed, assume $c(g) =0 \Rightarrow c(h) =0$, but there exists some cocycle, such that
$c(h) =0$ but $c(g) \neq 0$. If $c(g) \neq 0$, then $c(g) = (g-1)\xi$ for some $0 \neq \xi \in \cU G$ and the cocycle
$k \mapsto (k - 1)\xi - c(k) $ vanishes on $g$. This implies $c(h)=(h-1)\xi \neq 0$, which is a contradiction. 

If $c(g)=c(h)=0$, then also $c(gh)=0$ and $c(g^{-1})=0$, so that the equivalence classes together with the unit form subgroups. Denote by $[g]_{1} = \{h \in \dot G \mid h \sim g\} \cup \{e\}$. If $\beta_{1}^{(2)}([g]_{1}) \neq 0$, we continue
with the partitioning into subsets and proceed by transfinite induction. This implies claims $(i)$ and $(iv)$.

Claim $(ii)$ is proved by the following argument. 
If $hgh^{-1} \in [g]_{1}$ and $c(g)=0$, then $0=c(hgh^{-1}) = (1-hgh^{-1})c(h)$, and hence $c(h)=0$. We conclude
that $h \in [g]_{1}$.

We now prove $(iii)$ under the assumption of condition $(\star)$. Let $h=w_1v_1w_2v_2\cdots v_n$ be a shortest
(with respect to block length) trivial word consisting of non-trivial words $w_k \in G_i$ and $v_k \in G_j$ with $i \neq j$.
Let $c$ be a co-cycle which vanishes on $G_i$, but not on $G_j$. Then:
$$0 = c(h) = \left\{ w_1(v_1-1)+ \cdots + w_1v_1 \cdots w_n(v_n-1) \right\} \xi,$$
for some $\xi \in \cU G$.
We conclude by $(\star)$ that 
$$w_1(v_1-1)+ \cdots + w_1v_1 \cdots w_n(v_n-1) =0 \in \C G$$ and hence
there has to be a shorter trivial word. This contradicts the assumption and hence $G_i$ is free from $G_j$.
This finishes the proof.
\end{proof}

\begin{remark} \label{iinfinite}
It follows from Theorem \ref{bounded}, that the set $I$ is infinite if the first $\ell^2$-Betti number of $G$ does not
vanish. Indeed, if $I$ were finite, then $G$ would be boundedly generated by subgroups $G_i$ with vanishing first $\ell^2$-Betti number. This contradicts Theorem \ref{bounded}.
\end{remark}

The following lemma is useful to exploit the technique further.

\begin{lemma} \label{strong}
Let $G$ be a torsionfree discrete countable group satisfying condition $(\star)$. Let $H \subset G$ be a subgroup and assume that the restriction map
$${\res^G_H}\colon H^1(G,\cU G) \to H^1(H,\cU G)$$
is not injective. Then, there exists $h \in G$, such that the natural map
$$\pi\colon \Z \ast H \to \langle h,H \rangle \subset G$$
is an isomorphism. 
Moreover, if $G= \langle h,H \rangle$, then $h$ is free from $H$.
\end{lemma}
\begin{proof} Since ${\res^G_H}$ is not injective, there exists a non-trivial co-cycle $c \colon G \to \cU G$ which
is inner on $H$.
Subtracting this inner co-cycle, we can assume
that the restriction vanishes. Since $c$ was non-trivial, there exists $h \in G$, such that $c(h) = (h-1)\xi \neq 0$.
The proof proceeds as before.
\end{proof}

In the next two subsections we collect some corollaries of the results of the preceeding section.

\subsection{Simplicity of the reduced group $C^*$-algebra}

The following corollary shows that a non-trivial first Betti number of a torsionfree group
implies the existence of free subgroups. The proof uses the validity of condition $(\star)$ for the group.
It would be desirable to remove this assumption.

\begin{coro} \label{vonneumann}
Let $G$ be a discrete countable group satisfying condition $(\star)$. Assume that the first $\ell^2$-Betti number does not vanish.
Let $F$ be a finite subset of $G$. There exists $g \in G$, such that $g$ is free from each element in $F$. In particular, $G$ contains a copy of $F_2$.
\end{coro}
\begin{proof}
In view of the Theorem \ref{st1}, this can fail only if the index set $I$ in the proof of Theorem \ref{st1} is finite. Hence,
the result follows from Remark \ref{iinfinite}.
\end{proof}

\begin{remark} \label{insuff}
Corollary \ref{vonneumann} confirms the feeling that a sufficiently
non-amenable group contains a free subgroup. Note, that various weaker conditions like 
\emph{non-amenability} itself or \emph{uniform non-amenability} have been proved to be insufficient 
to ensure the existence of free subgroups, at least in the presence of torsion.
\end{remark}

Using results from \cite{MR1320606} we obtain the following result.

\begin{coro}
Let $G$ be a torsionfree discrete countable group satisfying condition $(\star)$. If the first $\ell^2$-Betti number
does not vanish, then the reduced group $C^*$-algebra $C^*_{red}(G)$ is simple and carries a unique trace.
\end{coro}
\begin{proof}
This follows from Lemma $2.2$ and Lemma $2.1$ 
in \cite{MR1320606}. Indeed, assuming condition $(\star)$ and non-vanishing
first $\ell^2$-Betti number, Corollary \ref{vonneumann} verifies Condition $P_{\rm nai}$ from 
Definition $4$ of \cite{MR1320606}.
\end{proof}

\subsection{Freiheitssatz and uniform exponential growth}

The following result is a generalization of the main result of J.\ Wilson in \cite{MR2011971} for torsionfree
groups which satisfy $(\star)$. For this, note that a group $G$ with $n$ generators and $m$ relations satisfies
$\beta_1^{(2)}(G) \geq n-m-1$.

\begin{coro}[Freiheitssatz] \label{manyfreesubgroups}
Let $G$ be a torsionfree discrete countable group which satisfies $(\star)$. 
Assume that $a_1,\dots,a_n \in G$ generate $G$ and $\lceil \beta^{(2)}_{1}(G) \rceil \geq k$. There exist $k+1$ elements
$a_{i_0}, \dots, a_{i_k}$ among the generators such that the natural map $$\pi \colon F_{k+1} \to \langle a_{i_0}, \dots, a_{i_k} \rangle \subset G$$ is an isomorphism.
\end{coro}
\begin{proof} We proof this result by induction over $n$. The case $n=1$ is obvious, since $n \geq k+1$ and there is
nothing to prove.
For the induction step, consider the restriction map
$$\res_1\colon H^1(G,\cU G) \to H^1\left(\langle a_2,\dots,a_n \rangle, \cU G\right).$$
If $\res_1$ is injective, then we can pass to the subgroup $G' = \langle a_2,\dots,a_n \rangle$. Note that $ \lceil \beta_1^{(2)}(G') \rceil \geq k$. In this case the proof is finished by induction since we decreased the number of generators by $1$. 

Hence, we can assume that the map $\res_1$ is not injective
and there exists a cocycle on $G$ which is inner on $G'=\langle a_2,\dots,a_n \rangle$. Lemma \ref{strong} implies
that $G= \langle a_1 \rangle \ast G'$. Now, the number of generators of $G'$ is $n-1$ and $\lceil \beta_1^{(2)}(G') \rceil \geq k-1$. Again, the proof is finished by induction.
\end{proof}

Following the work of J.\ Wilson in \cite{MR2011971}, this gives also an easy proof of Conjecture $5.14$ of M.\ Gromov in
\cite{MR1699320}, saying that the exponential growth rate of a group with $n$ generators and $m$ relations is bigger than
$2(n-m)-1$. Recall, the exponential growth rate is defined as
$$e_S(G) = \lim_{n \to \infty} \sqrt[n]{  \# B_S(e,n) },$$
where $B_S(e,n)$ denotes the ball of radius $n$ with respect to the word length metric coming from
a generating set $S$.
In general, we obtain the following result about the exponential growth rate:

\begin{coro}
Let $G$ be a finitely generated torsionfree discrete countable group which satisfies $(\star)$.
Then $$e_S(G) \geq 2 \lceil \beta_1^{(2)}(G) \rceil +1,$$
for any generating set $S$.
\end{coro}
\begin{proof} This is obvious, since Corollary \ref{manyfreesubgroups} says that a 
generating set $S$ contains the base of a free group of rank
$\lceil\beta_1^{(2)}(G) \rceil+1$. \end{proof}

In particular, a torsionfree group satisfying condition $(\star)$ has uniform exponential growth if its first
$\ell^2$-Betti number is positive.

\section{Results about the subgroup structure} \label{proofsec}

\subsection{Various notions of normality} \label{notionsofnormality}

We first want to review some notions of normality of subgroups and introduce some
notation. A subgroup $H \subset G$ is called:
\begin{enumerate}
\item[(i)] normal iff $gHg^{-1} =H$, for all $g \in G$,
\item[(ii)] $s$-normal iff $gHg^{-1} \cap H$ is infinite for all $g \in G$, and
\item[(iii)] $q$-normal iff $gHg^{-1} \cap H$ is infinite for elements $g \in G$, which generate $G$.
\end{enumerate}

We say that a subgroup inclusion $H \subset G$ satisfies one of the normality properties from above
\emph{weakly}, iff there exists an ordinal number $\alpha$, and an ascending $\alpha$-chain of subgroups, such that
$H_0 =H$, $H_{\alpha} = G$, and $\cup_{\beta < \gamma} H_\beta \subset H_{\gamma}$ has the required normality property.

Clearly, normal implies $s$-normal implies $q$-normal; and similarly for the weak notions.  The notion of (\emph{weakly}) $q$-normal 
subgroups were introduced by Popa (Definition 2.3 in \cite{MR2225044}) in order to ``untwist'' cocycles from the subgroup to the whole group
under certain weak mixingness conditions.
This method was also used quite successfully by Popa in subsequent works \cite{MR2231961, MR2231962, MR2342637}.
Theorem \ref{main1} below gives another instance where this notion is useful in untwisting cocycles.
See also Definition 1.2 in \cite{ioanapetersonpopa} for a von Neumann analogue of this notion.

A weakly $q$-normal 
subgroup is called $wq$-normal in \cite{MR2231962} and we follow this convention. In analogy, we call weakly $s$-normal subgroups
$ws$-normal. For obvious reasons $s$- and $q$-normality are considered only for infinite subgroups.

Weakly normal subgroups are usually called \emph{descendent}. Every subgroup is an descendent subgroup of a self-normalizing subgroup. Remark \ref{remark1} will clarify the corresponding observation in case of $wq$-normality. There are various
other notions of normality. For example,
P.\ Kropholler studies the notion of \emph{near normality} in \cite{kropp}. A subgroup $H \subset G$ is said to be near normal, if $H^g \cap H$ has finite index in $H$, for all $g \in G$. Clearly, near normality implies $s$-normality.

\begin{example}
The inclusions $$GL_n(\Z) \subset GL_n(\Q),\quad \mbox{and} \quad \Z = \langle x \rangle \subset 
\langle x,y \mid yx^p y^{-1} = x^q \rangle = BS_{p,q}$$ are inclusions of $s$-normal subgroups.
The inclusion $$F_2 = \langle a,b^2 \rangle \subset \langle a,b \rangle = F_2$$ is
$q$-normal but not $s$-normal.
\end{example}

\begin{lemma}[\cite{MR2225044}] \label{alter} Let $G$ be a discrete countable group and $H$ be an infinite subgroup.
The subgroup $H$ is $wq$-normal in $G$ if and only if given any intermediate subgroup
$H \subset K \subsetneq G$ there exists $g \in G \setminus K$ with
$gKg^{-1} \cap K$ infinite.
\end{lemma}
\begin{proof} One direction is obvious, since one can perform a transfinite induction to produce the chain of
subgroups with the desired properties.

We prove the converse:
Consider the least $\beta$, such that
$H_{\beta}$ is not contained in $K$. Then $\cup_{\gamma < \beta} H_{\gamma} \subset K$ and
there exists $g \in H_{\beta} \setminus K \subset G \setminus K$, such that $g\left(\cup_{\gamma < \beta} H_{\gamma}\right)g^{-1} \cap (\cup_{\gamma < \beta} H_{\gamma})$ is infinite. Hence,
$gKg^{-1} \cap K$ is infinite as well.
\end{proof}

\begin{remark} \label{remark1}
The notion of $wq$-normal subgroup is rather general. The following fact is easily deduced from the
previous lemma. If $G$ is torsionfree and
$H \subset G$, then $H$ is $wq$-normal in a malnormal subgroup of $G$. For general $G$, almost malnormal
subgroups have to be considered.
\end{remark}

\begin{coro} \label{second}
If $H\subset G$ is $wq$-normal, and $H \subset K \subset G$, then $K \subset G$ is $wq$-normal.
In particular if $H$ contains an infinite group which is normal in $G$ then $H$ is $wq$-normal.
\end{coro}
\begin{proof}
This is an immediate consequence of Lemma \ref{alter}
\end{proof}

We see from the next lemma, that $s$-normality shares slightly better inheritance properties than $q$-normality.
However, the following lemma does not seem to extend to the notion of $ws$-normality.

\begin{lemma} \label{lemma s-normal}
If $H \subset G$ is $s$-normal, and $H \subset K \subset G$, then $H \subset K$ and 
$K \subset G$ are inclusions of $s$-normal subgroups.
\end{lemma}
\begin{proof} This is obvious. \end{proof}

\subsection{$\ell^2$-invariants and normal subgroups}

The two main results in this subsection are Theorem \ref{main1} and Theorem \ref{main2}. We derive several
corollaries about the structure of groups $G$ with $\beta_1^{(2)}(G) \neq 0$.

\begin{theorem} \label{main1}
Let $G$ be a countable discrete group and suppose $H$ is an infinite $wq$-normal subgroup.
We have $\beta_1^{(2)}(H) \geq \beta_1^{(2)}(G)$.
\end{theorem}
\begin{proof} According to Theorem \ref{cohom}, the $\ell^2$-Betti-numbers are the $\cU G$-dimension of the
spaces $H^1(H,\cU G)$ and $H^1(G,\cU G)$. In order to show the inequality,
we show that the restriction map $H^1(G,\cU G) \to H^1(H,\cU G)$ is injective.
Let $c\colon G \to \cU G$ be a $1$-cocycle which is inner on $H$.
We may subtract the inner cocycle and assume that $c(h) = 0$, $\forall h \in H$.
Let $K = \{ g \in G | c(g) = 0 \}$, then $H \subset K \subset G$ and so if
$K \not= G$ then there exists $g \in G \setminus K$ with
$gKg^{-1} \cap K$ is infinite.
However, for each $k \in gKg^{-1} \cap K$
we have $c(g) - kc(g) = c(k) - gc(g^{-1}kg) = 0$. Hence, if
$g K g^{-1} \cap K$ is infinite we conclude that $c(g) = 0$. Indeed, this follows for $c(g) \in \ell^2 G$ from
strong mixing of the regular representation. The result extends to the general case by approximation in
rank metric. Thus $g \in K$ which gives a contradiction. Thus we conclude that $K = G$ which finishes the proof.
\end{proof}

\begin{remark} \label{remark2}
The inequality in Theorem \ref{main1} is sharp. Indeed $\langle a,b^2 \rangle \subset \langle a,b \rangle = F_2$ is $wq$-normal
and the restriction map in $\ell^2$-cohomology is an isomorphism. 
\end{remark}

\begin{coro} \label{coro2}
Let $H \subset K \subset G$ be a chain of subgroups and assume that $H \subset G$ is $wq$-normal and
$[K:H]< \infty$. Then
$$[K:H] \cdot  \beta_1^{(2)}(G) \leq \beta_1^{(2)}(H).$$
\end{coro}
\begin{proof}
This follows immediately from the proof of Theorem \ref{main1}, since the restriction map factorizes through
$H^1(K,\cU G)$. This $\cU G$-module has dimension $[K:H]^{-1} \beta_1^{(2)}(H)$, if the index $[K:H]$ is finite.
Alternatively, one can use Corollary \ref{second}.
\end{proof}

\begin{coro} \label{group}
Let $G$ be a torsionfree discrete countable group and let $H \subset G$ be an infinite subgroup. If $\beta_1^{(2)}(H)< \beta_1^{(2)}(G)$,
then there exists a proper malnormal subgroup $K \subset G$, such that $H \subset K$.
\end{coro}
\begin{proof}
This follows from Theorem \ref{main1} and Lemma \ref{alter}.
\end{proof}

\begin{remark}
Assume that $G$ is finitely presented of deficiency $d$ and that $H$ is finitely generated
with $n$ generators. Note that the hypothesis of Corollary \ref{group} is satisfied whenever $n < d$
or $H$ amenable and $0 < d$.
Moreover, the example in Remark \ref{remark2} shows that the assumption of a strict inequality 
cannot be improved.
\end{remark}

\begin{coro}
Let $G$ be a countable discrete group and let $H \subset G$ be an infinite $wq$-normal subgroup. Let 
$K \subset G$ be a subgroup with $H \subset K$ and assume that $\beta_1^{(2)}(G) > n$. 
Then, $K$ is not generated by $n$ or less elements.
\end{coro}

The second main result in this section is the following.

\begin{theorem} \label{main2}
Let $G$ be a countable discrete group and suppose $H$ is an infinite index, infinite $ws$-normal subgroup. If 
$\beta_1^{(2)}(H) < \infty$, then $\beta_1^{(2)}(G) =0$.
\end{theorem}
Obviously, for the proof we can restrict to the case of a $s$-normal subgroup. 
The proof of this theorem requires the introduction of some tools from ergodic theory and dynamical systems. It will be carried out in the next section. Note that the result follows from Theorem \ref{main1}, in case there are finite index subgroups
$G'$ of $G$, which have arbitrary high index and contain $H$. Indeed, in this case
$$ \beta^{(2)}_1(H) \geq \beta^{(2)}_1(G') = [G,G'] \cdot \beta^{(2)}_1(G),$$
by Theorem \ref{main1}, 
since Lemma \ref{lemma s-normal} implies that $H$ is also $s$-normal and hence $q$-normal in $G'$. 
This implies $\beta^{(2)}_1(G)=0$, under the assumption $\beta^{(2)}_1(H)< \infty$.

Although the existence of such families of finite index subgroups seems to be rare, it can always be achieved in the setting of discrete measured groupoids, see Lemma \ref{lem:finiteindex}. After having established the analogue of
Theorem \ref{main1} for discrete measured groupoids, the proof of Theorem \ref{main2} follows as before.

\begin{coro}\label{coromoon}
Let $G$ be a countable discrete group with $\beta_1^{(2)}(G) >0$. Suppose that $H \subset G$ is an infinite, finitely generated $ws$-normal subgroup. Then $H$ has to be of finite index.
\end{coro}

Note that the result applies in case $H$ contains an infinite normal subgroup. Hence, this result is a generalization of the classical results by A.\ Karass and D.\ Solitar in \cite{MR0086813}, H.\ Griffiths in \cite{MR0153832},
and B. Baumslag in \cite{MR0199247}. A weaker statement with additional hypothesis was proved as Theorem $1(2)$
in \cite{bmv}.

\begin{coro}[Gaboriau]
Let $G$ be a group with an infinite index, infinite, normal subgroup $H$ with $\beta_1^{(2)}(H)< \infty$, then $\beta^{(2)}_1(G)=0$.
\end{coro}

\begin{remark}
A generalization of Gaboriau's result to higher $\ell^2$-Betti numbers was obtained by R.\ Sauer and
the second author in \cite{sauerthom}. There it was shown that for a normal subgroup $N \subset G$ with 
all $\beta^{(2)}_p(N)=0$, for $p < q$, and
$\beta^{(2)}_q(N)$ finite, it follows that
$\beta^{(2)}_p(G)=0$, for $p \leq q$. 
The proof uses a Hochschild-Serre spectral sequence for discrete measured groupoids.
For more results in this direction, see \cite{sauerthom}.
\end{remark}

\section{Discrete measured groupoids} \label{ergodic}

\subsection{Infinite index subgroups}
After the statement of Theorem \ref{main2}, we outlined a proof in the presence of a descending chain
of finite index subgroups. In this subsection, we prove that such a chain exists as soon we pass to a suitable setting
of discrete measured groupoids.

\begin{lemma} \label{lemaction}
Let $G$ be a countable discrete group and let $H \subset G$ be a subgroup of infinite index.
There exists a standard probability space $(X,\mu)$ and an ergodic m.p.\ action of $G$ on $X$, such that the restriction of
the action to $H$ has a continuum of ergodic components.
\end{lemma}
\begin{proof}
We set $X = \prod_{gH \in G/H} [0,1]$ and let $\mu$ be the product of the Lebesgue measure. Then $G \curvearrowright X$
is ergodic, since $G \curvearrowright G/H$ is transitive with one infinite orbit. Moreover, the restriction of the action to $H$ leaves
the first factor invariant, and hence the space of ergodic components with respect to the $H$-action is continuous.
\end{proof}

The following lemma is another instance, where the flexibility of measure spaces and discrete measured groupoids
allows for constructions which are not possible in the realm of groups.

\begin{lemma} \label{lem:finiteindex} Let $G$ be a countable discrete group and let $H \subset G$ be a subgroup of infinite index.
There exists a standard probability space $(X,\mu)$, on which $G$ acts by m.p.\ Borel isomorphisms, such that the translation groupoid $X \rtimes G$ has finite index subgroupoids of arbitrary index which contain $X \rtimes H$.
\end{lemma}
\begin{proof}
Consider the space $(X,\mu)$ obtained from Lemma \ref{lemaction}. Consider a partition $Y= \cup_{i=1}^n Y_i$
of the space of ergodic components with respect to the action $H$. Assume that $\mu(Y_i) = n^{-1}$, for all $1 \leq i \leq n$.
Consider the subgroupoid $\calk \subset \calg$, which consists of those elements in $\calg$, which preserve the partition
of $Y$. Lemma $3.7$ of \cite{sauerthom} implies that the index of $\calk$ in $\calg$ is $n$.
\end{proof}

\subsection{Notions of normality for groupoids}
In \cite{sauerthom}, following the work of \cite{fsz}, the notion of strong normality of subgroupoids has been identified to be the right notion if one wants to
construct quotient groupoids. For our purposes, a weaker notion of normality is of importance.

\begin{definition}
Let $(\calg,\mu)$ be a discrete measured groupoid and let $\calr \subset \calg$ be a subgroupoid. The subgroupoid $\calr$ is
said to be $s$-normal, if for every local section $\phi$ of $\calg$, the set $\phi \calr \phi^{-1} \cap \calr $ has infinite measure.
The notion of $ws$-normality is defined similarly.
\end{definition}

The following lemma is the analogue of Lemma \ref{lemma s-normal} for discrete measured groupoids.  The proof is straightforward and we leave it as an exercise.

\begin{lemma}\label{lemgpoidsnormal}
Let $(\calg,\mu)$ be a discrete measure groupoid, let $A \subset \calg^0$ have positive measure and let $\calh \subset \calk
\subset \calg$ be subgroupoids. If $\calh \subset \calg$ is a $s$-normal inclusion, then $\calh_A \subset \calg_A$, $\calh
\subset \calk$ and $\calk \subset \calg$ are $s$-normal inclusions as well.
\end{lemma}

In order to relate $s$-normality for groups to $s$-normality for groupoids, we need the following technical lemma.

\begin{lemma}\label{lemreturn}
Let $G$ be an infinite countable discrete group which acts by m.p.\ Borel automorphisms on a probability space $(X, \mu)$.   Let $A \subset X$ be a Borel subset such that $\mu(A) \not= 0$ then $\limsup_{g \in G} \mu ( A \cap gA ) > 0$.
\end{lemma}

\begin{proof}
Take $0 < \varepsilon < \mu(A)$ and suppose that $F = \{ g \in G | \mu( A \cap gA) \geq \varepsilon \}$ is finite.  Then for all
$g \in G \setminus F$ we have $\mu(A^c \cap gA) \geq \mu(A) - \varepsilon$.  Let $n \in \N$ be such that $$n(\mu(A) - \varepsilon) > \mu(A^c)$$ and let $d > 0$ be such that 
$${n \choose 2} d < n(\mu(A) - \varepsilon) - \mu(A^c).$$  

Let $F' = \{ g \in G | \mu( A \cap gA ) \geq d \}$.
If $F'$ is finite, then we may take 
$g_1, \ldots, g_n \in G \setminus F$ such that $$g_j \not\in \bigcup_{i < j} g_iF', \quad \mbox{for all }  1\leq j \leq n.$$
Then 
\begin{eqnarray*}
\sum_{1 \leq i < j \leq n} \mu(g_iA \cap g_jA) &\geq& \sum_{1 \leq i \leq n} \mu(g_iA \cap A^c) - \mu(\cup_i g_iA \cap A^c)\\
& \geq& n(\mu(A) - \varepsilon) - \mu(A^c) \\
&>&  {n \choose 2} d.
\end{eqnarray*} 
Hence there exists $i < j$ such that $\mu(A \cap g_i^{-1}g_jA) \geq
d$. This contradicts the fact that $g_j \not\in g_iF'$ and hence we must have that $F'$ is infinite, 
i.e.\ $\limsup_{g \in G} \mu ( A \cap gA ) \geq d > 0.$
\end{proof}

The next theorem shows that an $s$-normal inclusion of groups leads to an $s$-normal inclusion of translation
groupoids.

\begin{theorem}
Let $G$ be a countable discrete group and let $H \subset G$ be a $ws$-normal subgroup. Moreover, let $(X,\mu)$ be a standard probability space on which $G$ acts by m.p.\ Borel automorphisms. Then, $X \rtimes H \subset X \rtimes G$ is an inclusion of a $ws$-normal
subgroupoid.
\end{theorem}
\begin{proof}
It is enough to treat the case of an $s$-normal subgroup. Each local section of $X \rtimes G$ is a countable sum of sections of
the form $\chi_A g$, where $A \subset X$ is Borel of positive measure and $g \in G$. It is enough to prove the assertion for
local sections of the form $\chi_Ag\colon g^{-1}A \to A$. In which case $(\chi_Ag)(X \rtimes H)(\chi_Ag)^{-1} \cap
(X \rtimes H)$ contains
$$\left(g^{-1}A \times \{g\}\right) \cdot \left(X \times \{h\}\right) \cdot \left(A \times \{g^{-1}\}\right) =
\left((A \cap (gh^{-1}g^{-1})A) \times \{ghg^{-1}\} \right), $$ for each $h \in H \cap g^{-1} H g$.  Hence its measure is infinite
by applying Lemma \ref{lemreturn} to the infinite group $H \cap g^{-1} H g$.
\end{proof}

\subsection{$\ell^2$-invariants of discrete measured groupoids}

We define the first complete $\ell^2$-cohomology of $\calg$ to be
$$H^1(\calg, \cU(\calg,\mu)) = \underline{\Ext}^1_{\calr(\calg,\mu)} \left( L^{\infty}(\calg^0),\cU(\calg)\right),$$
where $\underline{\Ext}$ denotes the derived functor from the abelian category of $L^{\infty}(X)$-complete $\calr(\calg)$-modules. All these notions were explained in great detail in \cite{sauerthom}. What is important for our purposes is the following concrete and familiar description of $H^1(\calg, \cU(\calg,\mu))$ as a space cocycles modulo inner cocycles.

A $\calg$-cocycle with values in $\cU(\calg,\mu)$ is an assignment $c$
of an element in $\cU(\calg,\mu)$ to every local section, such that

\begin{enumerate}
\item $c(\phi) \in \ran(\phi)\cU(\calg,\mu)$,
\item $c$ is compatible with countable decompositions, and
\item $c(\phi \circ \psi) = \phi \cdot c(\psi) +  \ran(\phi \circ \psi) \cdot c(\phi)$.
\end{enumerate}
A $\calg$-cocycle with values in $\cU(\calg,\mu)$ is said to be \emph{inner}, if
there exists $\xi \in \cU(\calg,\mu)$, such that $c(\phi) = (\phi  - \ran(\phi)) \cdot \xi $, for all local sections
$\phi$.

Clearly, the vector space of $\calg$-cocycles forms a right module over $\cU(\calg,\mu)$. We denote
this module by $\Der(\calg,\cU(\calg,\mu))$.

\begin{propo}
Let $(\calg,\mu)$ be an infinite discrete measured groupoid. The following sequence of $\cU(\calg,\mu)$-modules is
exact
$$0 \to \cU(\calg,\mu) \to \Der(\calg,\cU(\calg,\mu)) \to H^1(\calg,\cU(\calg,\mu)) \to 0.$$
\end{propo}
\begin{proof} The proof follows the standard arguments in group cohomology, which are used to identify
the first cohomology with the space of co-cycles modulo inner co-cycles.
\end{proof}

\begin{lemma}
Let $G$ be a countable discrete group and let $(X,\mu)$ be a probability space, on which $G$ acts by m.p.\ Borel
automorphisms. Then,
$$\beta_1^{(2)}(X \rtimes G,\mu) = \beta_1^{(2)}(G).$$
\end{lemma}
\begin{proof}
D.\ Gaboriau found this result for free actions in his groundbreaking work 
\cite{gab2}. Later, again through a process of algebraization, it turned
out that freeness is not needed. A proof can be found in \cite{thomrank} or \cite{sauerthom}.
\end{proof}

\subsection{The analogue of Theorem \ref{main1} for groupoids}

We are now proving the analogue of Theorem \ref{main1} in the setup of discrete measured groupoids.
In view of the remarks after the statement of Theorem \ref{main2}, this finishes the proof of Theorem \ref{main2}.

\begin{theorem} Let $(\calg,\mu)$ be a discrete
measured groupoid and let $\calh$ be a $ws$-normal subgroupoid.
Then the restriction map
$$H^1\left(\calg,\cU(\calg,\mu)\right) \to H^1\left(\calh,\cU(\calg,\mu)\right)$$ is injective.
\end{theorem}

\begin{proof}
For simplicity we will restrict to the case when $\calh$ is a $s$-normal subgroupoid.  Let $c$ be a $\calg$-cocycle with values
in $\cU(\calg, \mu)$ and suppose that $c$ is inner when restricted to $\calh$. By subtracting an inner cocycle we may assume
that $c(\phi) = 0$ for all local sections of $\calh$.

Let $\psi$ be a local section for $\calg$, and let $p$ be the 
maximal projection in $L^\infty(\calg^0)$ such that $pc(\psi) = 0$, note that $p \geq 1 - \ran(\psi)$. Set
$\chi_A = 1-p$. If $\chi_A \not=0$, then by Lemma \ref{lemgpoidsnormal} 
the inclusion $\calh_A \subset \calg_A$ is also $s$-normal and thus the
set $(\chi_A\psi)^{-1} \calh_A (\chi_A \psi) \cap \calh_A$ has infinite measure.  Thus there exist local sections $\phi_n$ for
$\calh_A$ which have large support (i.e.\ $\liminf_{n \to \infty} \mu(\ran(\phi_n)) \not= 0$), converge weakly to $0$ as partial
isometries acting on $L^2(L(\calg, \mu))$ and such that $(\chi_A\psi)^{-1} \phi_n \chi_A\psi$ is again a local section for
$\calh_A$.  By the Banach-Alaoglu theorem we may take a subsequence and suppose that $\ran(\phi_n)$ converges weakly to a non-zero
positive element $x \leq \chi_A$.

Since $L^2(L(\calg, \mu))$ is rank dense in $\cU(\calg,\mu)$ we may use the normality of the action and the cocycle relation to show that for a $\| \cdot \|_2$-dense collection of vectors $\xi \in L^2(L(\calg, \mu))$ we have 
\begin{eqnarray*}\langle x c(\psi), \xi \rangle &=& \lim_{n \to \infty} \langle r(\phi_n) c(\psi), \xi \rangle \\
& =& \lim_{n \to \infty} \langle c(\psi(\chi_A\psi)^{-1} \phi_n \chi_A\psi), \xi \rangle \\
& =& \lim_{n \to \infty} \langle \phi_n c(\psi), \xi \rangle = 0,
\end{eqnarray*} 
hence $x c(\psi) = 0$ contradicting the maximality of $p$ because $0 \lneq x \leq \chi A$.  Thus we must have that $c = 0$ which
completes the proof.
\end{proof}

\section{Applications} \label{misc}

In this section we collect applications of Theorem \ref{main1} and Theorem \ref{main2}. In particular, we find upper bounds on the first
$\ell^2$-Betti number of a group 
in terms of the first $\ell^2$-Betti numbers of the constituents of the group. We reprove and generalize several
results about non-existence of certain infinite index subgroups in certain situations.

\subsection{Boundedly generated groups}

Let $G$ be a discrete group and $G_1,\dots,G_n$ be sub-groups. The group $G$ is said to be \emph{boundedly
generated} by the subgroups $G_1,\dots,G_n$, if there exists an integer $k \in \N$, such that every element
in $G$ is a product of less than $k$ elements from $G_1,\dots,G_n$.
The following theorem is based on an idea of A.\ Ioana.

\begin{propo} \label{bounded}
Let $G$ be a countable discrete group. If $G$ is boundedly 
generated by subgroups $G_1,\dots,G_n$, then the following relation
holds:
$$\beta^{(2)}_1(G) \leq \sum_{i=1}^ n \beta^{(2)}_1(G_i).$$
\end{propo}
\begin{proof}
It suffices to show that the restriction map
$$H^1(G,\cU G) \to \oplus_{i=1}^n H^1(G_i,\cU G)$$ is injective. Indeed,
$ \dim_{LG}H^1(G_i,\cU G) = \beta^{(2)}_1(G_i),$
and therefore:
$$\dim_{LG} \oplus_{i=1}^n H^1(G_i,\cU G) = \sum_{i=1}^n \beta^{(2)}_1(G_i).$$
Let $c\colon G \to \cU G$ be a cocycle which is in the kernel.
The cocycle $c$ is inner on $G_i$. Being inner, we find a projection $p_i \in RG$ with $\tau(p_i^\perp) \leq \epsilon/n$,
such that $cp_i\colon G_i \to \cU G$ will be uniformly bounded in the $2$-norm. We consider $p = \inf p_i$. Since $G$
is boundedly generated by the $G_1,\dots,G_n$, we conclude, using the cocycle relation, that
$cp$ is uniformly bounded in the $2$-norm on $G$. 
Hence it is inner and there exists $\xi_{\epsilon} \in \ell^2 G$ such that
$c(g)p = (g-1)\xi_{\epsilon}p$, for all $g\in G$. It follows that $\xi_{\epsilon}$ converges in rank metric to some vector $\xi \in \cU G$
and $c(g) = (g-1) \xi$, for all $g \in G$. This finishes the proof.
\end{proof}

This generalizes Proposition $5$ in \cite{ma}. A particular case of the preceding theorem is $G=SL_n(\Z)$, for $n \geq 3$,
which is boundedly generated by copies of $\Z$.

\subsection{Certain groups generated by a family of subgroups}

Let $G$ be a countable discrete group and let $\{G_\alpha \mid \alpha \in V \}$ be a 
family of subgroups. We define a graph $V_G=(V,E)$ with vertices $V$ and an edge between
$\alpha$ and $\beta$, if and only if the intersection $G_{\alpha} \cap G_{\beta}$ is infinite.

\begin{propo} \label{family}
Let $G$ be a group and let $\{G_\alpha \mid \alpha \in V \}$ be a family of sub-groups. 
Assume that
\begin{enumerate}
\item[(i)] $\cup_{\alpha \in V} G_\alpha$ generates $G$ as a group, and
\item[(ii)] the graph $V_G$ is connected.
\end{enumerate}
Then,
$$\beta^{(2)}_1(G) \leq \sum_{\alpha \in V} \beta^{(2)}_1(G_i).$$
\end{propo}
\begin{proof}
Following the ideas in the proof of Theorem \ref{bounded}, one can show that the restriction map
$$H^1(G,\cU G) \to \oplus_{\alpha \in V} H^1(G_\alpha,\cU G)$$
is injective.
\end{proof}

\subsection{Limit groups} \label{limitgroups}

The notion of \emph{limit groups} or \emph{fully residually free groups} was introduced by Z.\ Sela in \cite{MR1863735}.
A countable discrete group $\Gamma$ is said to be a limit group, if it is finitely generated and 
for every finite set $T \subset \Gamma$, there
exists a group homomorphisms $\phi\colon \Gamma \to F_2$, which is injective on $T$. It was shown by C.\ Champetier
and V.\ Guirardel in \cite{MR2151593}, that limit 
groups are precisely the limits of free groups in R.\ Grigorchuk's space of marked
groups. 
Moreover, M.\ Pichot showed in \cite{MR2250857}, that a semi-continuity 
property holds for $\ell^2$-Betti numbers
in the space of marked groups. In particular, as noted in \cite{MR2250857}:
$$\beta_1^{(2)}(\Gamma) \geq 1, \quad \mbox{for all non-abelian limits groups $\Gamma$.}$$

Hence, our results apply and in particular a generalization of 
Theorem $3.1$ and Corollary $3.4$ in \cite{MR2262789} follows from Theorems \ref{main1} and \ref{main2}.

Another implication of our results are Theorems $B$ and $C$ of I.\ Kapovich in \cite{MR1859278}, which were 
partially withdrawn in \cite{MR1938758}.
Indeed Theorem $B$ follows from Corollary \ref{coromoon}. For convenience, we restate 
Theorem $C$ of \cite{MR1859278} in the generality to which we can extend this result:

\begin{propo} \label{finiteindex}
Let $G$ be a countable discrete group and let $H,K \subset G$ be finitely generated
infinite subgroups. Assume that 
$[H: H \cap K]$ and $[K:H \cap K]$ are finite. If $\beta^{(2)}_{1}(\langle H,K \rangle) \neq 0$, then the inclusion
$H \cap K \subset \langle H,K \rangle$ has finite index.
\end{propo}
\begin{proof}
We show that $H \subset \langle H,K \rangle$ is a $s$-normal inclusion. 
Note that $(H \cap K) \cap (H \cap K)^g$ is of finite index in $H \cap K$ for all elements $g \in \langle H,K \rangle.$
Hence $H \cap H^g$ is infinite, for all $g \in \langle H,K \rangle$. Since $H$ is $s$-normal and finitely 
generated, $[\langle H:K \rangle,H]$ has to be finite, by Corollary \ref{coromoon}. This finishes the proof.
\end{proof}

\subsection{Power-absorbing subgroups}

The following definition has been studied in \cites{MR0086813, MR0153832, MR1955619}. We reprove
most of the results and generalize to groups with non-vanishing first $\ell^2$-Betti number.

\begin{definition} Let $G$ be a torsionfree discrete countable group. A subgroup $H \subset G$ is
called \emph{power-absorbing}, if for every $g \in G$, there exists $n \in \N$, such that
$g^n \in H$.
\end{definition}

It has been studied under which conditions a finitely generated normal power-absorbing
subgroup has to be of finite index.

\begin{propo}
Let $G$ be a torsionfree discrete countable group and $H \subset G$ be a power-absorbing finitely generated subgroup.
If $\beta_1^{(2)}(G) \neq 0$, then the subgroup $H$ has to be of finite index.
\end{propo}
\begin{proof}
Clearly, $H$ is $s$-normal in $G$. Hence, the claim follows from Corollary \ref{coromoon}.
\end{proof}

\subsection{Groups measure equivalent to free groups}
\label{measeq}
In \cite{gabfree}, D.\ Gaboriau investigates groups which are measure equivalent to free groups.
Examples of such groups include amenable groups, lattices in $SL(2, \mathbb R)$, elementarily free groups (\cite{elemfree}), and is stable under taking free products.  Gaboriau in particular shows that this class is stable under taking subgroups and that a group in this class
with vanishing first $\ell^2$-Betti number is amenable.

Recall that if $H$ is a subgroup of $G$ then the \emph{commensurator subgroup} ${\rm comm}_G(H)$ of $H$ in $G$ is the group 
of all $g \in G$ such that $H \cap H^g$ has finite index in $H$ and $H^g$.  From the above proof we can show the following:

\begin{propo}
Let $G$ be measure equivalent to a free group and $H$ a finitely generated subgroup of $G$ then either ${\rm comm}_G(H)$ is amenable or else $H$ has finite index in ${\rm comm}_G(H)$.
\end{propo}

We can also show that groups in this class have unique maximal amenable extensions:

\begin{propo}
Let $G$ be measure equivalent to a free group and $H$ an infinite amenable subgroup, then $H$ has a unique maximal amenable extension.
\end{propo}

\begin{proof}
Suppose that $H_1$, and $H_2$ are amenable subgroups which contain $H$, we must show that $\langle H_1, H_2 \rangle$ is also amenable.
As $H$ is infinite it follows from Theorem \ref{family} that $$\beta_1^{(2)}(\langle H_1, H_2 \rangle) = 0$$ and hence $\langle H_1, H_2 \rangle$ is indeed amenable.
\end{proof}

Note that the above theorems will also hold for the class of groups which admit a proper $\ell^2$-co-cycle.
Also note that the above theorem should have a von Neumann algebra analog.  Specifically, it should be the case that if $Q \subset L\mathbb F_2$ is a diffuse amenable von Neumann subalgebra of the free group factor, then $Q$ has a unique maximal amenable extension in $L\mathbb F_2$.  Some evidence for this appears in \cite{pete}, \cite{jung}, and in \cite{ozpopa}.

\end{document}